\documentclass[twocolumn]{autart}    
\usepackage[authoryear]{natbib}
\usepackage{amsmath, amsfonts, amssymb}

\usepackage{graphicx}    
\usepackage{subcaption}  
\usepackage{float}       

\usepackage{tikz}
\usetikzlibrary{arrows.meta} 
\usetikzlibrary{positioning} 
\usetikzlibrary{calc}        
\usetikzlibrary{fit}         
\usetikzlibrary{backgrounds} 

\usepackage[colorlinks=true, citecolor=blue, linkcolor=blue, urlcolor=blue, bookmarks=false]{hyperref}

\newtheorem{theorem}{Theorem}
\newtheorem{lemma}{Lemma}

\newtheorem{corollary}{Corollary}
\newtheorem{definition}{Definition}

\newtheorem{remark}{Remark}

\begin{document} 

\begin{frontmatter}

\title{Explicit Convergence Regions of PID-Damped Accelerated Gradient Methods in Nonconvex Optimization} 


\thanks[footnoteinfo]{Corresponding author.}

\author[HDU]{Ailun Jian}\ead{allenjian805@163.com}, 
\author[PolyU]{Xun Li}\ead{li.xun@polyu.edu.hk}, 
\author[HDU]{Weigang Sun}\ead{wgsun@hdu.edu.cn}, 
\author[ZUST]{Gaohang Yu\thanksref{footnoteinfo}}\ead{maghyu@163.com}

\address[HDU]{School of Sciences, Hangzhou Dianzi University, Hangzhou, 310027, Zhejiang, China} 
\address[PolyU]{Department of Applied Mathematics, The Hong Kong Polytechnic University, Kowloon, Hong Kong, China} 
\address[ZUST]{School of Science, Zhejiang University of Science and Technology, Hangzhou, 310023, Zhejiang, China}

\begin{abstract}                          
Momentum-based accelerated gradient methods are widely adopted to expedite convergence in nonconvex optimization, but are prone to overshooting and oscillatory behavior. A class of PID-damped accelerated gradient methods mitigates this issue by augmenting classical momentum methods with a discrete-time derivative damping term. However, the coupling among the step size, momentum, and derivative gain renders the explicit characterization of their convergence regions analytically intractable, leaving explicit theoretical convergence boundaries unexplored. In this paper, we model this class of algorithms as a third-order nonlinear feedback dynamical system and establish explicit three-dimensional convergence regions for the step size, momentum, and derivative gain via a robust control-theoretic analysis based on the Kalman-Yakubovich-Popov (KYP) lemma, formally guaranteeing linear convergence under the Regularity Condition. Furthermore, we reveal a strict geometric upper bound on the derivative gain dictated by the nonconvex curvature, beyond which over-damping severely contracts the feasible step-size region, providing a rigorous theoretical explanation for the overdamped stagnation phenomenon. Numerical experiments corroborate the theoretical boundaries and illustrate practical parameter selection guidelines for the derivative gain.
\end{abstract}
\begin{keyword}
	Nonconvex optimization; Accelerated gradient methods; PID control; Convergence analysis; Regularity Condition; Kalman-Yakubovich-Popov lemma.
\end{keyword}
\end{frontmatter}

\section{Introduction}
Accelerated gradient methods incorporating momentum have been widely adopted to expedite convergence in nonconvex optimization~\citep{nesterov1998introductory, polyak1964some, su2014differential, ma2018quasi}. By accumulating historical gradient information, these methods can effectively escape shallow local minima. However, the momentum accumulation mechanism inevitably induces overshooting and oscillatory behavior along the optimization trajectory~\citep{ogata1995discrete, gitman2019understanding}. From a control-theoretic perspective, the classical Heavy Ball (HB) method operates as a Proportional-Integral (PI) optimizer, and the Nesterov Accelerated Gradient (NAG) method admits a similar PI structure with a look-ahead correction~\citep{polyak1964some, shi2020rethinking}. To mitigate the oscillatory dynamics inherent in PI-type momentum, a natural extension is to incorporate a derivative (D) damping term, yielding a three-term PID-structured accelerated gradient method~\citep{wang2020pid, shi2020rethinking}. This derivative term, defined as the gradient difference, acts as a discrete-time damping mechanism to suppress trajectory oscillations~\citep{weng2022adapid, chen2024accelerated}.

Despite the empirical effectiveness of PID-structured accelerated methods, their theoretical convergence properties remain poorly understood. Beyond the step size and momentum coefficient, these algorithms require tuning an additional damping coefficient that is strongly coupled with the other two parameters, making hyperparameter selection highly non-trivial. While tuning guidelines for classical PID controllers are well-established~\citep{ziegler1942optimum, astrom2006advanced, zhou2024active}, and explicit parameter manifolds have been derived for stabilizing continuous-time nonlinear systems~\citep{zhao2017pid}, these results do not directly translate to discrete-time optimization algorithms, as these continuous-time frameworks typically rely on Lyapunov methods for smooth flows and do not account for the discretization errors and memory-expansion effects inherent in multi-step optimization. To analyze such algorithms rigorously, robust control theory has emerged as a principled framework that establishes an equivalence between algorithmic convergence and the asymptotic stability of a nonlinear feedback dynamical system~\citep{lessard2016analysis, hu2017control, xiong2020analytical}. Within this paradigm, \citet{lessard2016analysis} utilized Integral Quadratic Constraints (IQCs) and Linear Matrix Inequalities (LMIs) to numerically verify convergence. Building upon this foundation, \citet{xiong2020analytical} applied the Kalman-Yakubovich-Popov (KYP) lemma and a system relaxation strategy under the Regularity Condition~\citep{candes2015phase, chi2019nonconvex, chen2019gradient} to derive explicit analytical convergence regions for the HB and NAG methods. However, these results are restricted to two-term momentum methods and do not extend to the three-term PID-structured case.

While existing robust control techniques successfully resolve standard momentum methods, integrating a discrete-time derivative damping term introduces fundamental mathematical challenges. Specifically, tracking historical gradients elevates the state space into a third-order dynamical system, transforming previously manageable inequalities into a highly coupled, analytically intractable frequency-domain problem. To overcome this structural barrier, we introduce a novel system relaxation strategy that explicitly derives robust convergence guarantees and the exact three-dimensional parameter manifold under the Regularity Condition. The main contributions of this paper are summarized as follows:

1) We establish explicit three-dimensional convergence regions for the hyperparameter triplet $(\alpha, \beta, \gamma)$ of 
PID-damped accelerated gradient methods under the Regularity Condition, providing deterministic parameter selection criteria and formally guaranteeing linear convergence in nonconvex optimization.

2) We develop a generalized system relaxation strategy within a robust control-theoretic framework, which reduces the high-dimensional frequency-domain inequalities arising from the third-order state-space into analytically tractable geometric boundaries. This methodology extends the KYP-based analysis to accelerated algorithms with augmented state spaces arising from historical gradient tracking.

3) We characterize the dual physical effects of the derivative damping gain from the derived theoretical boundaries. While the derivative term suppresses trajectory oscillations, it possesses a strict geometric upper bound dictated by the curvature of the nonconvex landscape, beyond which over-damping severely contracts the feasible step-size region. This provides a rigorous theoretical explanation for the overdamped stagnation phenomenon observed in practice.

The remainder of this paper is organized as follows. Section~\ref{sec:framework} formulates the PID-damped accelerated gradient methods as a discrete-time nonlinear feedback system under the Regularity Condition. Section~\ref{sec:stability} develops the system relaxation strategy to satisfy the KYP prerequisites. Section~\ref{sec:fdi_analysis} derives the explicit three-dimensional convergence regions. Section~\ref{sec:exp} validates the theoretical boundaries through numerical experiments. Section~\ref{sec:conclusion} concludes the paper.
\section{Problem Formulation and Dynamical System Modeling} \label{sec:framework}
Consider the general unconstrained problem
\begin{equation}\label{eq:opt_problem}
	\min_{\boldsymbol{x} \in \mathbb{R}^d} f(\boldsymbol{x}),
\end{equation}
where $f:\mathbb{R}^d\to\mathbb{R}$ is a continuously differentiable nonconvex function. The analysis is conducted under the Regularity Condition (RC), formally defined as follows.
\begin{definition}[Regularity Condition] \label{def:RC}
	A differentiable function $f: \mathbb{R}^d \to \mathbb{R}$ satisfies the Regularity Condition $RC(\mu, \lambda, \epsilon)$ if there exist positive constants $\mu, \lambda$ and a neighborhood radius $\epsilon > 0$ such that, for a local minimizer $\boldsymbol{x}_*$, the inequality
	\begin{equation} \label{eq:RC_ineq}
		\langle \nabla f(\boldsymbol{x}), \boldsymbol{x} - \boldsymbol{x}_* \rangle \ge \frac{\mu}{2} \|\nabla f(\boldsymbol{x})\|^2 + \frac{\lambda}{2} \|\boldsymbol{x} - \boldsymbol{x}_*\|^2
	\end{equation}
	holds for all $\boldsymbol{x} \in \mathcal{N}_\epsilon(\boldsymbol{x}_*) := \{ \boldsymbol{x} : \|\boldsymbol{x} - \boldsymbol{x}_*\| \le \epsilon \|\boldsymbol{x}_*\| \}$. The Cauchy-Schwarz inequality restricts the parameters to $\mu\lambda \le 1$.
\end{definition}

\begin{remark}\label{rem:local_to_global}
	We initially assume the RC holds for all $\boldsymbol{x}\in\mathbb{R}^{d}$, denoted as $RC(\mu, \lambda)$. Lemma~\ref{lem:local_rc_init} establishes an explicit initialization condition on $\boldsymbol{x}_0$ and $\boldsymbol{x}_{-1}$ that confines all subsequent iterates within the neighborhood $\mathcal{N}_\epsilon(\boldsymbol{x}_*)$, validating the stability boundaries derived under this premise for the local regime.
\end{remark}

As noted by \citet{chi2019nonconvex}, the RC characterizes a combination of one-point strong convexity and smoothness. It ensures that the negative gradient direction consistently provides a descent path towards the minimizer $\boldsymbol{x}_*$ without requiring global convexity. Equipped with this algebraic boundary, we proceed to model PID-damped accelerated gradient methods as a third-order nonlinear feedback system.


\subsection{Unified Third-Order Feedback System Modeling} \label{sec:RC}
The discrete-time formulation of the classical PID controller is described by the recursive control sequence
\begin{equation} \label{eq:pid_position}
	\boldsymbol{u}_k = K_p \boldsymbol{e}_k + K_i \sum_{j=0}^{k} \boldsymbol{e}_j + K_d (\boldsymbol{e}_k - \boldsymbol{e}_{k-1}),
\end{equation}
where $\boldsymbol{u}_k$ is the control action and $\boldsymbol{e}_k$ represents the tracking error at step $k$. Applying the backward difference $\boldsymbol{u}_k-\boldsymbol{u}_{k-1}$ and the gradient mapping $\boldsymbol{e}_k = -\nabla f(\boldsymbol{x}_k)$ to \eqref{eq:pid_position} yields the PID optimization rule~\citep{shi2020rethinking}
\begin{equation}\label{eq:pid_update}
	\begin{split}
		\boldsymbol{x}_{k+1} &= \boldsymbol{x}_k - K_p \nabla f(\boldsymbol{x}_k) + K_i (\boldsymbol{x}_k - \boldsymbol{x}_{k-1}) \\
		&\quad - K_d (\nabla f(\boldsymbol{x}_k) - \nabla f(\boldsymbol{x}_{k-1})).
	\end{split}
\end{equation}
To generalize this algorithm into a unified PID-damped accelerated gradient framework, we introduce the parameter vector $\boldsymbol{\xi} = (\alpha, \beta, \hat{\beta}, \gamma)^\top$ to formulate the following discrete-time update
\begin{equation}\label{eq:unified_sys}
	\begin{aligned}
		\boldsymbol{x}_{k+1} &= (1+ \beta) \boldsymbol{x}_{k} - \beta \boldsymbol{x}_{k-1} - (\alpha + \gamma)\nabla f(\boldsymbol{y}_k) \\
		&\quad + \gamma \nabla f(\boldsymbol{y}_{k-1}),\\
		\boldsymbol{y}_k &= (1+\hat{\beta}) \boldsymbol{x}_{k} - \hat{\beta} \boldsymbol{x}_{k-1}.
	\end{aligned}
\end{equation}
Here, $\alpha > 0$ is the step size, $\beta \ge 0$ and $\hat{\beta} \ge 0$ denote the standard and look-ahead momentum coefficients respectively, and $\gamma \ge 0$ represents the damping coefficient of derivative term. By defining the nonlinear gradient feedback as $\boldsymbol{u}_k^{(1)} = \nabla f(\boldsymbol{y}_k)$, the unified recurrence \eqref{eq:unified_sys} can be recast as a discrete-time nonlinear feedback dynamical system
\begin{equation}
	\left\{
	\begin{aligned}
		\boldsymbol{x}_{k+1}^{(1)} &= (1 + \beta)\boldsymbol{x}_{k}^{(1)}- \beta \boldsymbol{x}_{k}^{(2)} - (\alpha + \gamma) \boldsymbol{u}_k^{(1)} + \gamma \boldsymbol{u}_{k}^{(2)}\\
		\boldsymbol{x}_{k+1}^{(2)} &= \boldsymbol{x}_{k}^{(1)}\\
		\boldsymbol{y}_k &= (1 + \hat{\beta}) \boldsymbol{x}_{k}^{(1)} - \hat{\beta} \boldsymbol{x}_{k}^{(2)}\\
		\boldsymbol{u}_{k}^{(1)} &= \nabla f(\boldsymbol{y}_k)\\
		\boldsymbol{u}_{k+1}^{(2)} &= \boldsymbol{u}_{k}^{(1)}
	\end{aligned}
	\right.
\end{equation}
where the variables are defined as $(\boldsymbol{x}_{k}^{(1)}, \boldsymbol{x}_{k}^{(2)}, \boldsymbol{u}_k^{(1)}, \boldsymbol{u}_k^{(2)}) := (\boldsymbol{x}_{k}, \boldsymbol{x}_{k-1}, \boldsymbol{u}_k, \boldsymbol{u}_{k-1})$. In this state-space configuration, $\boldsymbol{x}_{k}^{(2)}$ and $\boldsymbol{u}_k^{(2)}$ function as delay states, embedding historical position and gradient information into the current update rule. This unified parameterization recovers several standard accelerated methods via specific configurations $\boldsymbol{\xi}$
\begin{equation}\label{eq:order-accelerated method}
	\left\{
	\begin{aligned}
		& \text{GD}, & \boldsymbol{\xi} &= (\alpha, 0, 0, 0),\\
		& \text{HB}, & \boldsymbol{\xi} &= (\alpha, \beta, 0 , 0),\\
		& \text{NAG}, & \boldsymbol{\xi} &= (\alpha, \beta, \hat{\beta}, 0), \\
		& \text{PID}, & \boldsymbol{\xi} &= (\alpha, \beta, 0 , \gamma).
	\end{aligned}
	\right.
\end{equation}
Specifically, enforcing $\gamma = 0$ reduces our generalized framework to the classical momentum algorithms. As shown in Fig.~\ref{fig:feedback_system}, by applying the Kronecker product to decouple the spatial dimension $d$, the core linear time-invariant (LTI) system $G(\mathbf{A}, \mathbf{B}, \mathbf{C}, \mathbf{D})$ is formulated as
\begin{equation}\label{eq:LTI_model}
	\begin{aligned}
		\boldsymbol{\phi}_{k+1} &= (\mathbf{A}\otimes \mathbf{I}_d)\boldsymbol{\phi}_{k} + (\mathbf{B}\otimes \mathbf{I}_d)\boldsymbol{u}_k\\
		\boldsymbol{y}_k &= (\mathbf{C}\otimes \mathbf{I}_d)\boldsymbol{\phi}_{k} + (\mathbf{D}\otimes \mathbf{I}_d)\boldsymbol{u}_k.
	\end{aligned}
\end{equation}
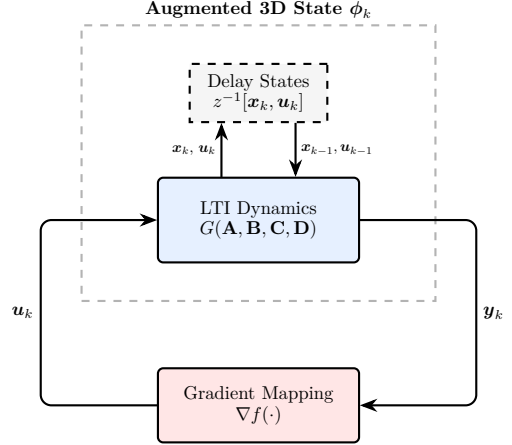
\begin{figure}[htbp] 
	\centering
	\begin{tikzpicture}[>=Stealth, thick, scale=0.7, every node/.style={transform shape}] 
	
		\definecolor{sysblue}{RGB}{230, 240, 255}
		\definecolor{sysred}{RGB}{255, 230, 230}
		\definecolor{sysgray}{RGB}{245, 245, 245}
		
		\node[draw, fill=sysblue, minimum width=3.8cm, minimum height=1.6cm, align=center, rounded corners=2pt] (G) at (0, 2) {LTI Dynamics \\ $G(\mathbf{A}, \mathbf{B}, \mathbf{C}, \mathbf{D})$};
		
		\node[draw, fill=sysred, minimum width=3.8cm, minimum height=1.4cm, align=center, rounded corners=2pt] (nabla) at (0, -1.5) {Gradient Mapping \\ $\nabla f(\cdot)$};
		
		\node[draw, dashed, fill=sysgray, minimum width=2.6cm, minimum height=1.1cm, align=center] (delay) at (0, 4.4) {Delay States\\ $z^{-1}[\boldsymbol{x}_{k}, \boldsymbol{u}_{k}]$};
		
		\draw[->] ([xshift=-0.7cm]G.north) -- ([xshift=-0.7cm]delay.south) node[midway, left, scale=0.8] {$\boldsymbol{x}_k$, $\boldsymbol{u}_k$};
		\draw[<-] ([xshift=0.7cm]G.north) -- ([xshift=0.7cm]delay.south) node[midway, right, scale=0.8] {$\boldsymbol{x}_{k-1}, \boldsymbol{u}_{k-1}$};
		
		\draw[->, rounded corners] (G.east) -- ++(2.2,0) |- (nabla.east)
		node[pos=0.25, right] {$\boldsymbol{y}_k$};
		
		\draw[->, rounded corners] (nabla.west) -- ++(-2.2,0) |- (G.west)
		node[pos=0.25, left] {$\boldsymbol{u}_k$};
		
		\begin{scope}[on background layer]
			\node[draw, dashed, gray!60, thick, inner xsep=1cm, inner ysep=0.5cm, fit=(G) (delay), 
			label={[anchor=south, inner sep=3pt, scale=1, font=\bfseries]north:Augmented 3D State $\boldsymbol{\phi}_k$}] (box) {};
		\end{scope}
	\end{tikzpicture}
	\caption{Dynamical system representation of PID-damped accelerated gradient methods. The augmented LTI system~\eqref{eq:LTI_model} incorporates a delay operator $z^{-1}$ to track historical states $(\boldsymbol{x}_{k-1}, \boldsymbol{u}_{k-1})$, forming a three-dimensional state-space formulation}
	\label{fig:feedback_system}
\end{figure}
The composite three-dimensional state $\boldsymbol{\phi}_{k}$ and the input $\boldsymbol{u}_k$ are constructed as
\begin{equation}
	\boldsymbol{\phi}_{k} = [{\boldsymbol{x}_{k}^{(1)}}^\top, {\boldsymbol{x}_{k}^{(2)}}^\top, {\boldsymbol{u}_k^{(2)}}^\top]^\top, \quad
	\boldsymbol{u}_k = \boldsymbol{u}_k^{(1)}.
\end{equation}
The state-space matrices governing the internal dynamics are given by
\begin{equation}\label{eq:matrix}
	\left[\begin{array}{c|c}
		\mathbf{A} & \mathbf{B} \\
		\hline 
		\mathbf{C} & \mathbf{D}
	\end{array}\right]
	= \left[\begin{array}{ccc|c}
		1+\beta & -\beta & \gamma & -(\alpha+\gamma) \\
		1 & 0 & 0 & 0 \\
		0 & 0 & 0 & 1 \\
		\hline 
		1+\hat{\beta} & -\hat{\beta} & 0 & 0
	\end{array}\right].
\end{equation}

Let $\boldsymbol{\phi}^{*} = [{\boldsymbol{x}^{*}}^\top, {\boldsymbol{x}^{*}}^\top, \mathbf{0}^\top]^\top \in \mathbb{R}^{3d}$ denote the equilibrium point of \eqref{eq:LTI_model}. This state corresponds strictly to a stationary point $\boldsymbol{x}^{*}$ of the function $f$, where the gradient vanishes. Since the asymptotic stability of the feedback system ($\boldsymbol{\phi}_{k} \to \boldsymbol{\phi}^{*}$) guarantees the algorithmic convergence ($\boldsymbol{x}_{k} \to \boldsymbol{x}^{*}$), verifying optimization convergence is equivalent to analyzing the stability of this dynamical system.

\subsection{Linear Convergence via the KYP Lemma} \label{sec:lyapunov_stability}

Building upon the LTI state-space model~\eqref{eq:LTI_model} established in Section~\ref{sec:framework}, we leverage the RC~\eqref{eq:RC_ineq} to formulate a quadratic constraint.

\begin{lemma}[Quadratic Constraint]\label{lem:quad_constraint}
	Suppose the objective function $f$ satisfies the $\text{RC}(\mu, \lambda)$ with a stationary point $\boldsymbol{x}_{*}$. By defining the output error $\tilde{\boldsymbol{y}}_k = \boldsymbol{y}_k - \boldsymbol{y}_*$ and the gradient error $\tilde{\boldsymbol{u}}_k = \boldsymbol{u}_k - \boldsymbol{u}_*$, the regularity condition on $f$ guarantees that the nonlinear feedback mapping satisfies the following quadratic constraint
	\begin{equation}\label{eq:wrc}
		\begin{bmatrix}
			\tilde{\boldsymbol{y}}_k \\
			\tilde{\boldsymbol{u}}_k
		\end{bmatrix}^\top
		\left(\mathbf{M} \otimes \mathbf{I}_d\right)
		\begin{bmatrix}
			\tilde{\boldsymbol{y}}_k \\
			\tilde{\boldsymbol{u}}_k
		\end{bmatrix}\geq 0,
	\end{equation}
	where the symmetric multiplier matrix $\mathbf{M} \in \mathbb{R}^{2 \times 2}$ is exactly defined as
	\[
	\mathbf{M} = 
	\begin{bmatrix}
		-\lambda & 1 \\
		1 & -\mu
	\end{bmatrix}.
	\]
\end{lemma}

By integrating constraint of inequality~\eqref{eq:wrc} with the KYP framework, we can establish the linear convergence of the system through a time-domain LMI.

\begin{lemma}[Linear Convergence via LMI]\label{lem:stability}
	Suppose the objective function $f$ satisfies $\text{RC}(\mu, \lambda)$. Consider the system matrices $\mathbf{A}, \mathbf{B}, \mathbf{C}, \mathbf{D}$ defined in \eqref{eq:matrix}. If there exists a positive definite matrix $\mathbf{P}\succ 0$ and a scalar decay rate $\rho\in(0,1)$ satisfying the following LMI
	\begin{equation}\label{eq:LMI}
		\begin{bmatrix}
			\mathbf{C} & \mathbf{D} \\
			\mathbf{0} & 1
		\end{bmatrix}^\top
		\mathbf{M}
		\begin{bmatrix}
			\mathbf{C} & \mathbf{D} \\
			\mathbf{0} & 1
		\end{bmatrix}
		+
		\begin{bmatrix}
			\mathbf{A}^\top \mathbf{P} \mathbf{A}-\rho^2 \mathbf{P} & \mathbf{A}^\top \mathbf{P} \mathbf{B} \\
			\mathbf{B}^\top \mathbf{P} \mathbf{A} & \mathbf{B}^\top \mathbf{P} \mathbf{B}
		\end{bmatrix}\preceq 0,
	\end{equation}
	where $\mathbf{0}$ is a zero vector of appropriate dimensions, then the system state $\boldsymbol{\phi}_{k}$ converges linearly to the equilibrium point $\boldsymbol{\phi}_{*}$ with rate $\rho$
	\begin{equation*}
		\|\boldsymbol{\phi}_{k} - \boldsymbol{\phi}_{*}\|\leq \sqrt{\kappa(\mathbf{P})}\rho^k\|\boldsymbol{\phi}_{0} - \boldsymbol{\phi}_{*}\|,
	\end{equation*}
	where $\kappa(\mathbf{P}) = \lambda_{\max}(\mathbf{P}) / \lambda_{\min}(\mathbf{P})$ denotes the condition number of $\mathbf{P}$.
\end{lemma}

\begin{remark}\label{rem:kyp_motivation}
	Lemma~\ref{lem:stability} reduces the convergence analysis to finding a feasible $3 \times 3$ positive definite matrix $\mathbf{P}$. However, solving the resulting $4 \times 4$ LMI~\eqref{eq:LMI} algebraically to extract explicit parameter boundaries is intractable. To bypass this matrix inequality, we leverage the KYP lemma~\citep{rantzer1996kalman}, which transforms the LMI into a tractable scalar Frequency-Domain Inequality (FDI).
\end{remark}

\begin{lemma}[KYP Lemma]\label{lem:kyp}
	Given matrices $\mathbf{A}\in\mathbb{R}^{n\times n}$, $\mathbf{B}\in\mathbb{R}^{n\times m}$, and a symmetric matrix $\mathbf{\Pi}\in\mathbb{R}^{(n+m)\times (n+m)}$, assume that $\det(e^{j\omega}\mathbf{I} - \mathbf{A})\neq0$ for all $\omega\in\mathbb{R}$. The following two statements are equivalent:
	\begin{enumerate}
		\item The FDI holds for all $\omega \in \mathbb{R}$:
		\begin{equation}\label{eq:FDI}
			\begin{bmatrix}
				(e^{j \omega} \mathbf{I}-\mathbf{A})^{-1} \mathbf{B} \\
				\mathbf{I}
			\end{bmatrix}^* \mathbf{\Pi}
			\begin{bmatrix}
				(e^{j \omega} \mathbf{I}-\mathbf{A})^{-1} \mathbf{B} \\
				\mathbf{I}
			\end{bmatrix} \prec 0,
		\end{equation}
		where $*$ denotes the conjugate transpose.
		\item There exists a symmetric matrix $\mathbf{P}\in\mathbb{R}^{n\times n}$ such that the following LMI holds:
		\begin{equation}\label{eq:LMI_generic}
			\mathbf{\Pi}+
			\begin{bmatrix}
				\mathbf{A}^\top \mathbf{P} \mathbf{A}-\mathbf{P} & \mathbf{A}^\top \mathbf{P} \mathbf{B} \\
				\mathbf{B}^\top \mathbf{P} \mathbf{A} & \mathbf{B}^\top \mathbf{P} \mathbf{B}
			\end{bmatrix} \prec 0.
		\end{equation}
	\end{enumerate}
\end{lemma}

While Lemma~\ref{lem:kyp} guarantees the existence of a symmetric matrix $\mathbf{P}$, the convergence proof in Lemma~\ref{lem:stability} requires $\mathbf{P}$ to be strictly positive definite. To guarantee this positive definiteness, \citet{xiong2020analytical} established the exact prerequisites in the following corollary.

\begin{corollary}[KYPC Corollary]\label{coro:kypc}
	The strict FDI condition \eqref{eq:FDI} is equivalent to the existence of a positive definite matrix $\mathbf{P} \succ 0$ satisfying the LMI \eqref{eq:LMI_generic}, provided that the pair $(\mathbf{A}, \mathbf{\Pi})$ satisfies a set of prerequisites, denoted as $\text{KYPC}(\mathbf{A}, \mathbf{\Pi})$. Specifically, $\text{KYPC}(\mathbf{A}, \mathbf{\Pi})$ requires that
	\begin{enumerate}
		\item $\det(e^{j\omega}\mathbf{I} - \mathbf{A})\neq0$ for all $\omega\in\mathbb{R}$;
		\item The matrix $\mathbf{A}$ is strictly Schur stable;
		\item The upper-left block of the multiplier $\mathbf{\Pi}$ is positive semidefinite.
	\end{enumerate}
\end{corollary}

Corollary~\ref{coro:kypc} reduces the stability analysis to verifying $\text{KYPC}(\mathbf{A},\mathbf{\Pi})$ and evaluating the scalar FDI~\eqref{eq:FDI}. Next, we introduce the relaxation strategies required for the system $G(\mathbf{A}, \mathbf{B}, \mathbf{C}, \mathbf{D})$ to satisfy these prerequisites.

\section{System Relaxation for KYP Feasibility} \label{sec:stability}

The open-loop state matrix $\mathbf{A}$ formulated in~\eqref{eq:matrix} has an eigenvalue on the unit circle, violating the strict Schur stability prerequisite of the KYP lemma. To address this, we introduce a relaxation parameter $\delta < 0$ to shift the poles inward, reformulating the original dynamics into a strictly Schur stable system.

\subsection{Relaxation of PID-Damped Dynamics}

We define a relaxed feedback signal $\boldsymbol{v}_k = -\gamma\nabla f(\boldsymbol{y}_k) - \delta \boldsymbol{x}_k$ to incorporate the relaxation parameter $\delta$. Substituting the inverse relation $\nabla f(\boldsymbol{y}_k) = -\frac{1}{\gamma}(\boldsymbol{v}_k + \delta \boldsymbol{x}_k)$ into the update equation~\eqref{eq:unified_sys} yields the relaxed open-loop system
\begin{equation}
	\left\{
	\begin{aligned}
		\boldsymbol{x}_{k+1} &= \left(1 + \beta + \frac{(\alpha+\gamma)\delta}{\gamma}\right)\boldsymbol{x}_{k} - (\beta+\delta) \boldsymbol{x}_{k-1} \\
		&\quad - \boldsymbol{v}_{k-1} + \frac{\alpha + \gamma}{\gamma} \boldsymbol{v}_k, \\
		\boldsymbol{y}_k &= (1 + \hat{\beta}) \boldsymbol{x}_{k} - \hat{\beta} \boldsymbol{x}_{k-1}.
	\end{aligned}
	\right.
\end{equation}

To avoid increasing the system dimension, we set $\hat{\beta} = 0$, yielding $\boldsymbol{y}_k = \boldsymbol{x}_k$. However, this does not imply that the system excludes the NAG method. Through a linear substitution into~\eqref{eq:unified_sys} under $\gamma=0$, we obtain
\begin{equation} \label{eq:nag_expansion}
	\begin{aligned}
		\boldsymbol{y}_{k+1} &= \boldsymbol{x}_{k+1} + \beta(\boldsymbol{x}_{k+1} - \boldsymbol{x}_k) \\
		&= \underbrace{\boldsymbol{y}_k - \alpha\nabla f(\boldsymbol{y}_k) + \beta (\boldsymbol{y}_k - \boldsymbol{y}_{k-1})}_{\text{HB Momentum}} \\
		&\quad - \underbrace{\alpha\beta \big( \nabla f(\boldsymbol{y}_k) - \nabla f(\boldsymbol{y}_{k-1}) \big)}_{\text{Gradient Difference (D) Term}}.
	\end{aligned}
\end{equation}
This indicates that under a linear transformation, NAG can be viewed as an HB method augmented with a derivative term, rendering it a specific subclass derived from our system. Applying the identity $\boldsymbol{y}_k \equiv \boldsymbol{x}_k$, the augmented state vector is updated to $\boldsymbol{\phi}_k = [{\boldsymbol{x}_k}^\top, {\boldsymbol{x}_{k-1}}^\top, {\boldsymbol{v}_{k-1}}^\top]^\top$. Consequently, the relaxation reformulates the state update into the LTI system $G(\mathbf{\hat{A}}, \mathbf{\hat{B}}, \mathbf{\hat{C}}, \mathbf{\hat{D}})$ with input $\boldsymbol{v}_k$ and output $\boldsymbol{y}_k$. The corresponding block matrices are given by
\begin{equation}\label{eq:relaxed_matrices}
	\left[\begin{array}{c|c}
		\mathbf{\hat{A}} & \mathbf{\hat{B}} \\
		\hline 
		\mathbf{\hat{C}} & \mathbf{\hat{D}}
	\end{array}\right]
	=
	\left[\begin{array}{ccc|c}
		1 + \beta + \frac{(\alpha+\gamma)\delta}{\gamma} & - (\beta+\delta) & -1 & \frac{\alpha + \gamma}{\gamma} \\
		1 & 0 & 0 & 0 \\
		0 & 0 & 0 & 1 \\
		\hline 
		1 & 0 & 0 & 0 
	\end{array}\right].
\end{equation}
Since the relaxation alters the nonlinear feedback mapping, a modified quadratic constraint is required for the new input-output pairs.

\begin{lemma}[Modified Quadratic Constraint]\label{lem:relaxed_qc}
	Suppose $f$ satisfies the $\text{RC}(\mu, \lambda)$ with a stationary point $\boldsymbol{y}_{*}$. Under the feedback $\boldsymbol{v}_k = -\gamma\nabla f(\boldsymbol{y}_k)-\delta \boldsymbol{y}_{k}$, the input-output error pairs $\tilde{\boldsymbol{y}}_k = \boldsymbol{y}_k - \boldsymbol{y}_*$ and $\tilde{\boldsymbol{v}}_k = \boldsymbol{v}_k - \boldsymbol{v}_*$ satisfy the modified quadratic constraint
	\begin{equation}\label{eq:relaxed_qc}
		\begin{bmatrix}
			\tilde{\boldsymbol{y}}_k \\
			\tilde{\boldsymbol{v}}_k
		\end{bmatrix}^\top
		\left(\mathbf{\hat{M}} \otimes \mathbf{I}_d\right)
		\begin{bmatrix}
			\tilde{\boldsymbol{y}}_k \\
			\tilde{\boldsymbol{v}}_k
		\end{bmatrix} \geq 0,
	\end{equation}
	where the symmetric multiplier matrix $\mathbf{\hat{M}} \in \mathbb{R}^{2 \times 2}$ is defined as
	\[
	\mathbf{\hat{M}} = \begin{bmatrix}
		r_1  & r_2 \\
		r_2 & -\mu 
	\end{bmatrix},
	\]
	with $r_1 = -(\gamma^2\lambda + 2\gamma\delta +\mu\delta^2)$ and $r_2 = - (\gamma + \mu\delta)$.
\end{lemma}

\subsection{Feasibility of the Relaxed System}
We establish the exact conditions under which the relaxed system satisfies the KYPC.
\begin{theorem}[KYP Feasibility Conditions]\label{thm:param_range}
	Let $f$ satisfy the $\text{RC}(\mu, \lambda)$. There exists a relaxation parameter $\delta$ such that the relaxed system~\eqref{eq:relaxed_matrices} satisfies $\text{KYPC}(\mathbf{\hat{A}}, \mathbf{\hat{M}})$ if and only if the algorithmic parameters $\alpha$, $\beta$, and $\gamma$ satisfy
	\begin{equation}\label{eq:param_condition}
		\left\{
		\begin{aligned}
			&0 < \gamma < \frac{(1+\beta)(1+\sqrt{1-\mu\lambda})}{\lambda},\\
			&0 < \alpha < \frac{2(1+\beta)(1+\sqrt{1-\mu\lambda})}{\lambda}- 2\gamma,\\
			&-1<\beta<1,
		\end{aligned}
		\right.
	\end{equation}
	and the relaxation parameter $\delta$ is chosen from the interval
	\begin{equation}\label{eq:delta}
		\delta \in \left(\max\left\{-1-\beta, -\frac{2\gamma(1+\beta)}{2\gamma+\alpha}\right\}, 0\right).
	\end{equation}
\end{theorem}

\begin{remark}
	Theorem~\ref{thm:param_range} concludes the global deterministic stability analysis. It demonstrates that any algorithmic parameter set $\boldsymbol{\xi}$ within the derived boundaries guarantees a non-empty feasible interval for the relaxation parameter $\delta$. The existence of this relaxation parameter ensures the relaxed system satisfies the KYPC. Combining this result with Corollary~\ref{coro:kypc} and the LMI in Lemma~\ref{lem:stability} guarantees that the accelerated gradient method achieves linear convergence to the local minimizer $\boldsymbol{x}_*$.
\end{remark}

Theorem~\ref{thm:param_range} assumes the $\text{RC}(\mu,\lambda)$ without spatial restriction. Since general non-convex functions typically satisfy the RC only locally around a minimizer $\boldsymbol{x}_*$, the following lemma introduces an exact initialization condition to confine the trajectory within this valid neighborhood.

\begin{lemma}[Local RC Initialization]\label{lem:local_rc_init}
	Suppose a local minimizer $\boldsymbol{x}_*$ of $f$ satisfies the local $\text{RC}(\mu, \lambda, \epsilon)$. Let $\mathbf{P} \succ 0$ be a feasible solution to the LMI established in Lemma~\ref{lem:stability}. Define the neighborhood $\mathcal{B}_r(\boldsymbol{x}_*) := \{\boldsymbol{x} : \|\boldsymbol{x} - \boldsymbol{x}_*\| \le r\}$. If the initial points $\boldsymbol{x}_{-1}$ and $\boldsymbol{x}_0$ are chosen within the ball $\mathcal{B}_{r_0}(\boldsymbol{x}_*)$ with radius
	\begin{equation}
		r_0 = \frac{\epsilon \|\boldsymbol{x}_*\|}{\sqrt{\kappa(\mathbf{P})\left(2 + \left(\frac{2\gamma}{\mu} + |\delta|\right)^2\right)}},
	\end{equation}
	then the output sequence satisfies $\boldsymbol{y}_k \in \mathcal{B}_{\epsilon\|\boldsymbol{x}_*\|}(\boldsymbol{x}_*)$ for all $k \ge 0$.
\end{lemma}

Consequently, Lemma~\ref{lem:local_rc_init} validates the application of our global stability framework to local non-convex settings, ensuring the algorithm rigorously achieves local linear convergence.

\section{Explicit Convergence Regions} \label{sec:fdi_analysis}

\subsection{Frequency-Domain Inequality Transformation}\label{subsec:fdi_derivation}

Evaluating the scalar FDI~\eqref{eq:FDI} analytically requires the explicit formulation of the open-loop transfer function.

\begin{lemma}[Open-Loop Transfer Function]\label{lem:transfer_function}
	For the relaxed LTI system $G(\mathbf{\hat{A}}, \mathbf{\hat{B}}, \mathbf{\hat{C}}, \mathbf{\hat{D}})$ defined in~\eqref{eq:relaxed_matrices}, the open-loop transfer function $H(z) = \mathbf{\hat{C}}(z \mathbf{I}-\mathbf{\hat{A}})^{-1} \mathbf{\hat{B}} + \mathbf{\hat{D}}$ is given by
	\begin{equation}\label{eq:transfer_func_explicit}
		H(z) = \frac{B_1 z - 1}{z^2 - A_{11}z - A_{12}},
	\end{equation}
	where the intermediate scalars are defined as
	\begin{equation}\label{eq:sys_scalars}
		A_{11} = 1 + \beta + \frac{(\alpha+\gamma)\delta}{\gamma}, \quad
		A_{12} = -(\beta+\delta), \quad
		B_1 = \frac{\alpha + \gamma}{\gamma}.
	\end{equation}
\end{lemma}
\begin{theorem}\label{thm:fdi_polynomial}
	Let $f$ satisfy the $\text{RC}(\mu, \lambda)$, and suppose the parameters $\boldsymbol{\xi}$ satisfy Theorem~\ref{thm:param_range} to guarantee the KYPC. The existence of a matrix $\mathbf{P} \succ 0$ solving the LMI in Lemma~\ref{lem:stability} is equivalent to the quadratic polynomial inequality
	\begin{equation}\label{eq:polynomial_P_x}
		P(x) = \Gamma_2 x^2 + \Gamma_1 x + \Gamma_0 < 0
	\end{equation}
	holding for all $x = \cos\omega \in [-1, 1]$, where the coefficients are defined as
	\begin{equation}\label{eq:gamma_coeffs}
		\left\{
		\begin{aligned}
			\Gamma_2 &= 4(\gamma - \mu\beta), \\
			\Gamma_1 &= 2\lambda\gamma(\alpha+\gamma) - 2(1+\beta)(\alpha + 2\gamma) + 2\mu(1+\beta)^2, \\
			\Gamma_0 &= -\lambda(\alpha^2 + 2\alpha\gamma + 2\gamma^2) + 2\alpha(1+\beta) \\
			&\quad + 4\gamma\beta - 2\mu(1+\beta^2).
		\end{aligned}
		\right.
	\end{equation}
\end{theorem}

Theorem~\ref{thm:fdi_polynomial} reduces the high-dimensional LMI to a scalar quadratic constraint over a compact interval. Based on this equivalence, solving $P(x) < 0$ explicitly yields the exact three-dimensional analytical boundaries $(\alpha, \beta, \gamma)$ for the PID-damped accelerated gradient algorithm.

\subsection{Three-Dimensional Convergence Boundaries}
\begin{theorem}\label{thm:analytical_convergence_region}
	Let $f$ satisfy the $\text{RC}(\mu, \lambda)$. The PID-damped accelerated gradient algorithm achieves linear convergence to the local minimizer $\boldsymbol{x}_*$ if the algorithmic parameters $(\alpha, \beta, \gamma)$ belong to the feasible region $\Omega = \Omega_u \cup \Omega_d$, comprising
	
	1) The over-damped region $\Omega_u$, defined as
	\begin{equation}\label{eq:region_u}
		\Omega_u = \left\{ (\alpha, \beta, \gamma) \;\middle|\;
		\begin{aligned}
			&\beta \in (-1, 1), \\
			&\max\{0, \mu\beta\} \le \gamma < \frac{\mu(1+\beta)}{1+\sqrt{1-\mu\lambda}}, \\
			&0 < \alpha <U(\beta,\gamma),
		\end{aligned}
		\right\}.
	\end{equation}
	
	2) The under-damped region $\Omega_d$, defined as
	\begin{equation}\label{eq:region_d}
		\Omega_d = \left\{ (\alpha, \beta, \gamma) \;\middle|\;
		\begin{aligned}
			&\beta \in (0, 1), \\
			&0 < \gamma < \mu\beta, \\
			&0 < \alpha < \min \left\{ U(\beta,\gamma), \; \Xi_d(\beta, \gamma) \right\},
		\end{aligned}
		\right\},
	\end{equation}
	where $U(\beta,\gamma)= \frac{2(1+\beta)(1-\sqrt{1-\mu\lambda})}{\lambda} - 2\gamma$, and the implicit boundary is $\Xi_d(\beta, \gamma) := \{ \alpha > 0 \mid \Delta(\alpha) < 0 \}$, with $\Delta(\alpha) := \Gamma_1^2 - 4\Gamma_2\Gamma_0$.
\end{theorem}

\begin{remark}\label{rem:physical_interpretation}
	Theorem~\ref{thm:analytical_convergence_region} explicitly characterizes the stability domain of the PID-damped method, providing a theoretical justification for its empirical advantages in overshoot suppression.
	
	(i) It establishes deterministic selection criteria for the parameter triple $(\alpha, \beta, \gamma)$ strictly based on the objective function's geometric properties $(\mu, \lambda)$, explicitly revealing their algebraic coupling.
	
	(ii) For a fixed momentum $\beta \in (0, 1)$, the derivative gain $\gamma$ effectively mitigates the overshoot inherited from the HB momentum. However, the upper bound $U(\beta, \gamma)$ reveals a strict linear trade-off: increasing $\gamma$ necessitates a proportional reduction in the maximum allowable step size $\alpha$.
	
	(iii) The partition of $\Omega$ into $\Omega_u$ and $\Omega_d$ reflects distinct dynamic behaviors. The region $\Omega_u$ guarantees non-oscillatory over-damped convergence, whereas $\Omega_d$ accommodates low-derivative configurations by enforcing the stricter, implicit step-size bound $\Xi_d$ to ensure under-damped stability.
	
	(iv) As $\gamma \to 0^+$, the system degenerates into the classical HB method. Consequently, the over-damped region $\Omega_u$ vanishes, and the step-size upper bound exactly recovers the analytical limit established in~\citep{xiong2020analytical}. This follows directly from the definition of $U(\beta, \gamma)$ by taking $\gamma \to 0^+$, yielding $\lim_{\gamma \to 0^+} U(\beta, \gamma) = \frac{2(1+\beta)(1-\sqrt{1-\mu\lambda})}{\lambda}$.
\end{remark}
\section{Numerical Experiments} \label{sec:exp}
\begin{figure*}[t!]
	\centering
	\includegraphics[width=0.75\linewidth]{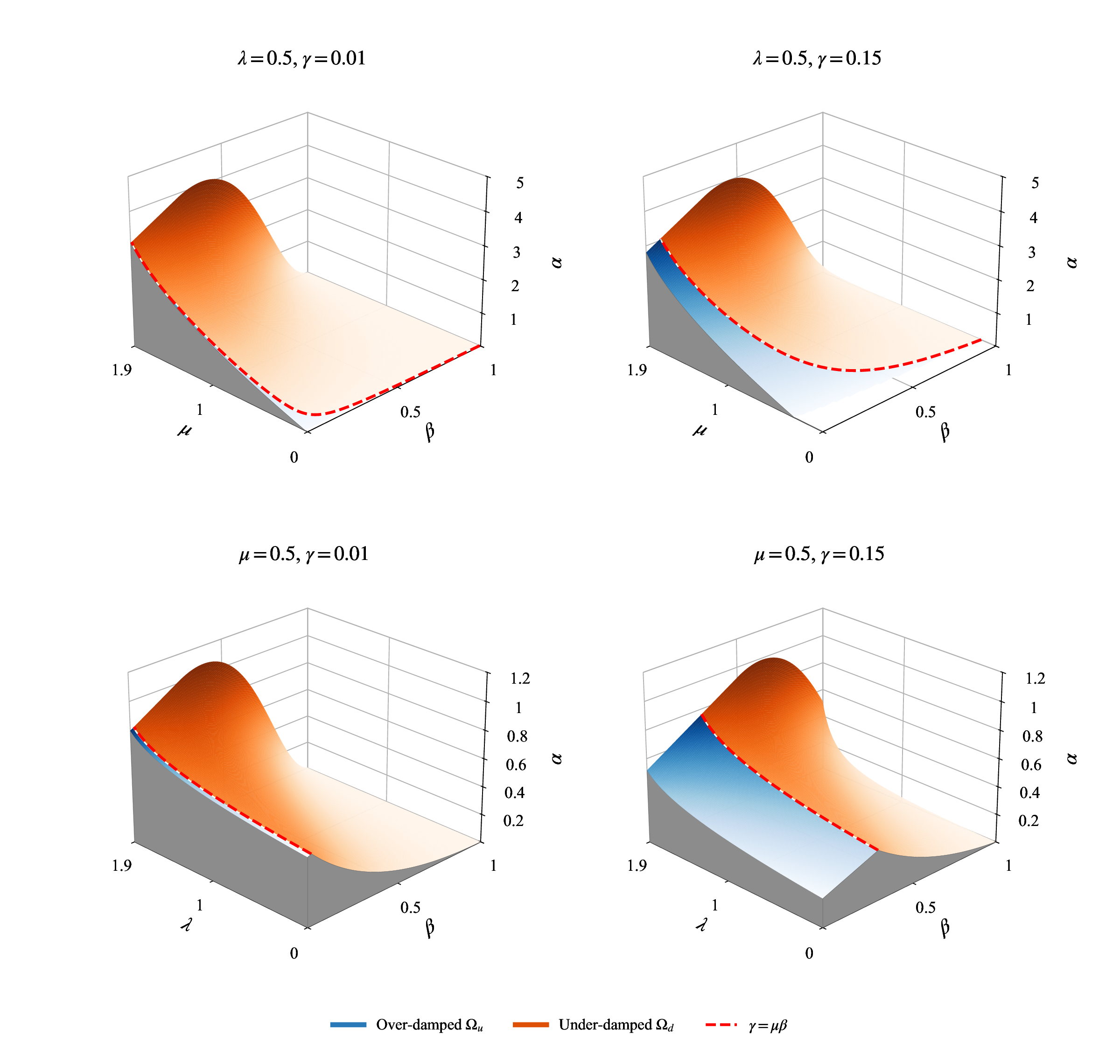}
	\caption{3D Convergence Manifolds of the PID-damped method. The $2 \times 2$ subplots illustrate the maximum allowable step size $\alpha$ under varying landscape parameters ($\mu$ and $\lambda$). The blue and orange surfaces represent the over-damped region $\Omega_u$ and the under-damped region $\Omega_d$, respectively, separated by the boundary $\gamma = \mu\beta$ (red dashed line).}
	\label{fig:3d_hb_pid}
\end{figure*}
In this section, we visualize the three-dimensional parameter domain established in Theorem~\ref{thm:analytical_convergence_region} and subsequently demonstrate the effects of the under-damped and over-damped regions on overshoot mitigation through a simple numerical example.

\subsection{Visualizations of Three-Dimensional Convergence Manifolds}

Based on the analytical boundaries defined in~\eqref{eq:region_u} and~\eqref{eq:region_d}, we determine the feasible parameter space for the PID-damped accelerated gradient algorithm and visualize it as a three-dimensional manifold. This explicitly illustrates the geometric trade-offs among $\alpha$, $\beta$, and $\gamma$ dictated by the theoretical boundaries.

Fig.~\ref{fig:3d_hb_pid} visualizes the three-dimensional convergence manifolds of the maximum allowable step size $\alpha$ under varying landscape parameters ($\mu$ and $\lambda$). For a negligible derivative gain ($\gamma=0.01$), the manifold reduces to the classical HB parameter domain, strictly dominated by the under-damped region $\Omega_d$. Increasing $\gamma$ to $0.15$ introduces a distinct over-damped region $\Omega_u$ that effectively mitigates overshoot. However, this reveals a fundamental algorithmic trade-off: to maintain stability, the introduction of $\gamma$ inevitably compresses the feasible space for $\alpha$, particularly under high-$\beta$ and low-$\mu$ configurations.

\begin{figure*}[!t]
	\centering
	\includegraphics[width=0.75\linewidth]{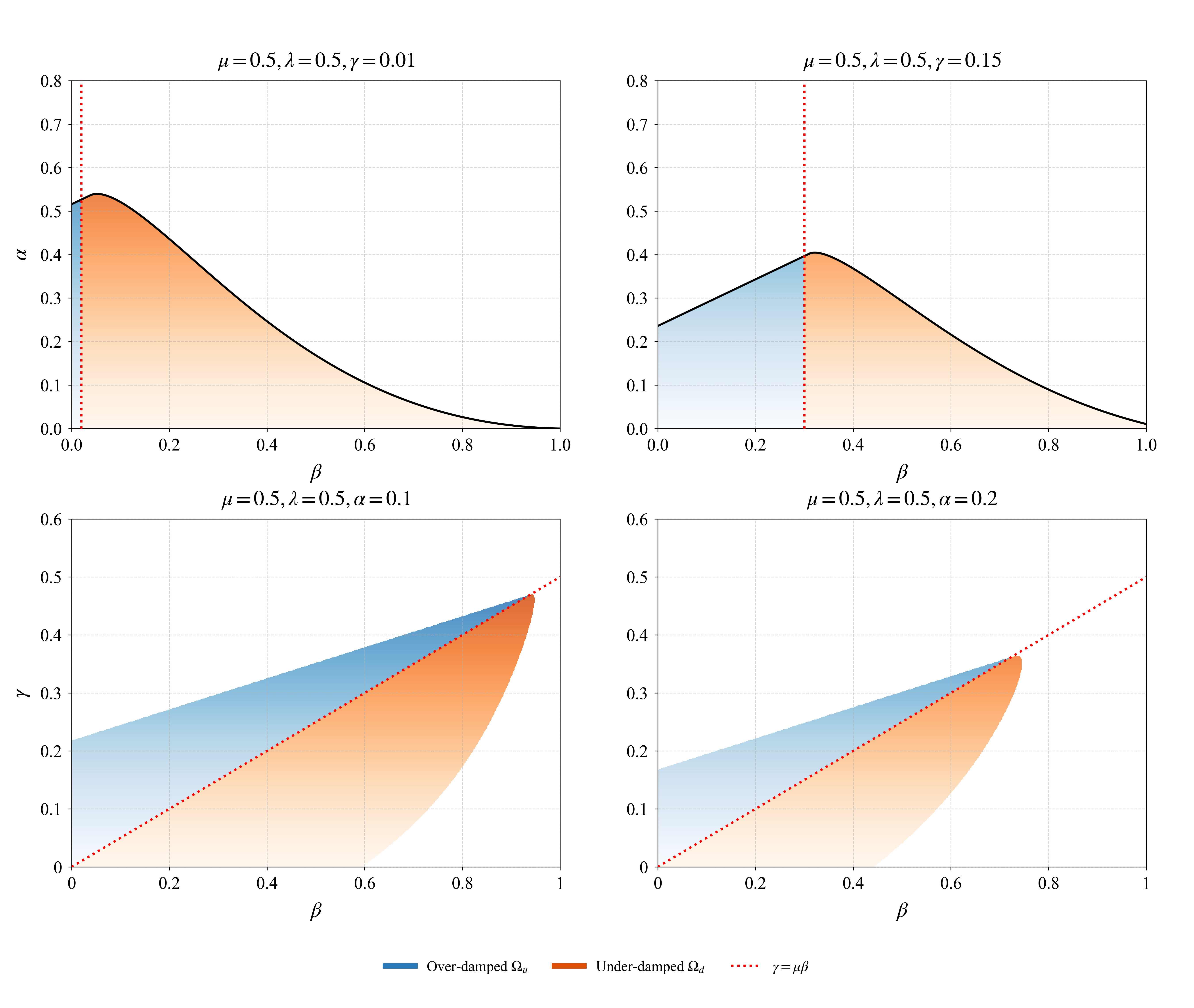}
	\caption{2D Parameter Cross-sections and Phase Diagrams of the PID-damped method ($\mu=0.5, \lambda=0.5$). The $2 \times 2$ subplots illustrate the allowable step size $\alpha$ versus momentum $\beta$ under fixed $\gamma$, and the derivative gain $\gamma$ versus momentum $\beta$ under fixed $\alpha$. The blue and orange regions represent the over-damped region $\Omega_u$ and the under-damped region $\Omega_d$, respectively, separated by the boundary $\gamma = \mu\beta$ (red dotted line).}
	\label{fig:2d_analysis}
\end{figure*}

To further decouple the parameter interactions, Fig.~\ref{fig:2d_analysis} provides 2D cross-sections and hyperparameter phase diagrams under a fixed landscape geometry ($\mu=0.5, \lambda=0.5$). The $\alpha$-$\beta$ cross-sections visually confirm that while a larger derivative gain $\gamma$ expands the over-damped proportion, it simultaneously imposes a linear penalty on the peak of the allowable step size. Furthermore, the $\gamma$-$\beta$ phase diagrams explicitly outline the safe operational corridors for practical design. They demonstrate that a more aggressive step size (e.g., $\alpha=0.2$) severely compresses the feasible parameter space, verifying that high momentum strictly requires a proportionally high derivative gain to prevent system divergence.



\subsection{Validation on an Ill-Conditioned Quadratic Valley}
To empirically validate the topological phase transitions and the overshoot mitigation mechanism, we evaluate the algorithms on an ill-conditioned quadratic valley. The objective function is defined as
\begin{equation}\label{eq:test_func}
	f(\boldsymbol{x}) = 2x_1^2 + 5x_2^2.
\end{equation}
For this landscape, the eigenvalues of the Hessian are $\lambda_{\min}	=4$ and $\lambda_{\max}=10$. According to the sector bound equivalence~\citep{lessard2016analysis,xiong2020analytical}, the geometric parameter is rigorously determined as $\mu = \frac{2}{L+m} = \frac{1}{7}$. Given a fixed momentum $\beta = 0.9$, the theoretical critical damping boundary for the phase transition is calculated as $\gamma^* = \mu\beta \approx 0.129$. The optimization processes are uniformly initialized at $\boldsymbol{x}_0 = [-2.0, -2.0]^\top$ and executed for 500 iterations. We benchmark the PID-damped method against GD, HB, and NAG. For all accelerated algorithms, the step size and momentum are fixed at $\alpha=0.01$ and $\beta=0.9$, respectively. To explore the damping effects, we evaluate the derivative gain across a discrete spectrum $\gamma \in \{0.01, 0.129, 0.30, 0.45, 0.70\}$.

\begin{figure*}[!t]
	\centering
	\includegraphics[width=0.75\linewidth]{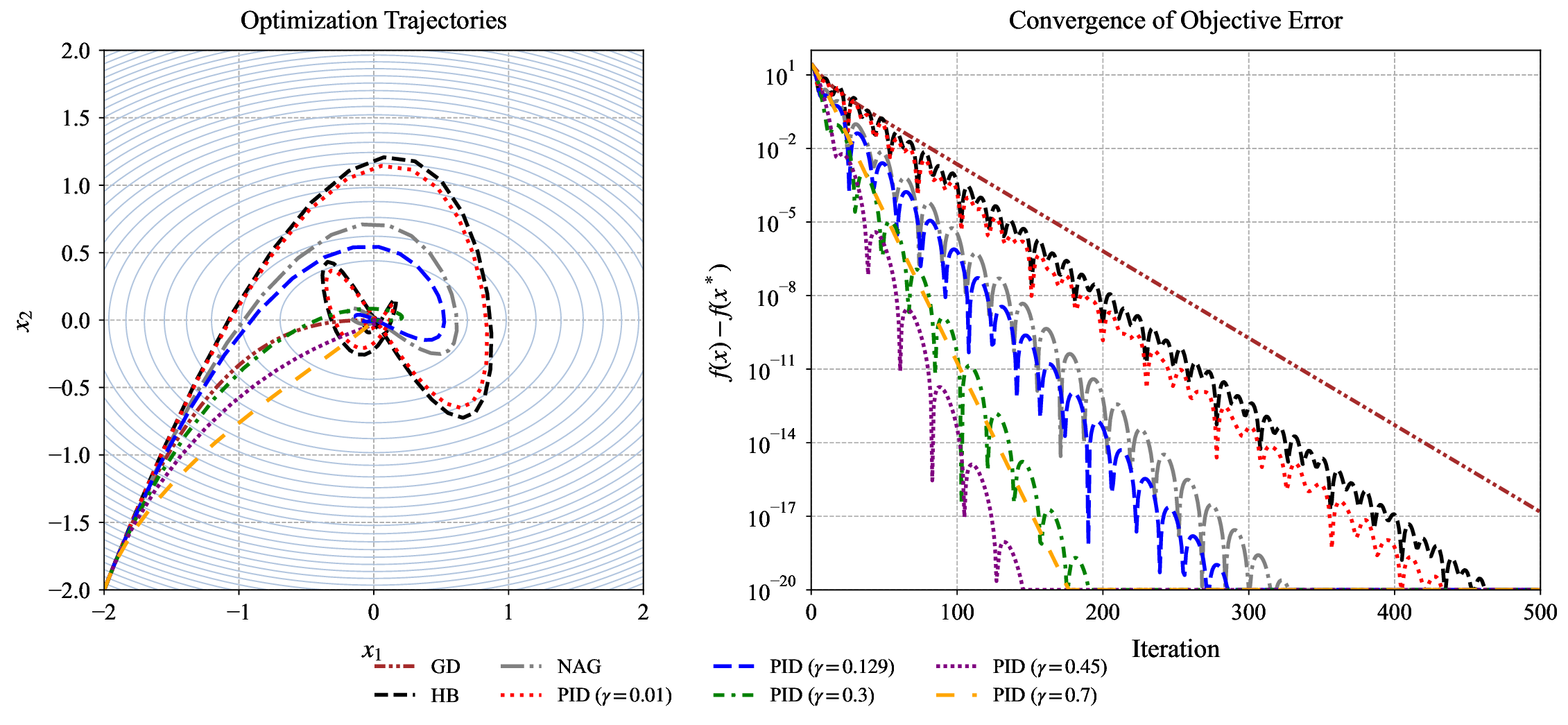}
	\caption{Empirical validation of topological phase transitions on an ill-conditioned quadratic valley. The $1 \times 2$ subplots illustrate the 2D optimization trajectories (left) and the objective error convergence (right) under varying derivative gains $\gamma$.}
	\label{fig:test_function}
\end{figure*}

Fig.~\ref{fig:test_function} illustrates the 2D optimization trajectories and the corresponding objective error convergence, which perfectly corroborate the theoretical boundaries. Specifically, when operating below the critical boundary ($\gamma = 0.01 < \gamma^*$), the system resides in the under-damped region $\Omega_d$ and exhibits pronounced Lissajous oscillations akin to the classical HB method. As $\gamma$ crosses the theoretical boundary $\gamma^* \approx 0.129$ and enters the over-damped region $\Omega_u$ (e.g., $\gamma \ge 0.15$), the severe oscillations are immediately suppressed, validating the topological phase transition. Within this stable $\Omega_u$ region, increasing $\gamma$ initially straightens the trajectory and accelerates convergence, identifying $\gamma=0.45$ as yielding the fastest empirical convergence within the over-damped region. However, an excessively large gain ($\gamma=0.70$) introduces an over-damping penalty, which severely suppresses the momentum effect and reduces the algorithm to near-gradient-descent behavior.

This demonstrates that crossing the theoretical boundary $\gamma^*$ is strictly necessary to eliminate overshoot, but finding the optimal convergence rate requires careful tuning within the over-damped corridor.

\section{Conclusion}\label{sec:conclusion}

This paper establishes explicit three-dimensional convergence regions for PID-damped accelerated gradient methods under the Regularity Condition, resolving an open problem in the convergence analysis of three-term accelerated algorithms. By modeling the algorithm as a third-order nonlinear feedback dynamical system and introducing a novel system relaxation strategy, the high-dimensional LMI feasibility problem is reduced to a scalar frequency-domain inequality, from which the analytical parameter boundaries are extracted in closed form. A key theoretical finding is the identification of a strict geometric upper bound on the derivative gain: while the derivative term effectively suppresses momentum-induced oscillations, exceeding this bound triggers overdamped stagnation by severely contracting the feasible step-size region. This provides the first rigorous explanation for a phenomenon widely observed but previously uncharacterized in derivative-damped optimization.

Several directions remain open for future investigation. First, the current analysis establishes the feasible convergence region but does not characterize the optimal parameter configuration within it; deriving the convergence rate as a function of $(\alpha, \beta, \gamma)$ and identifying the rate-optimal triplet within $\Omega$ is a natural next step. Second, extending the present framework to stochastic nonconvex optimization, where gradient observations are corrupted by noise, would broaden its applicability to large-scale machine learning settings. Third, the derived theoretical boundaries suggest a principled basis for adaptive damping schemes that dynamically adjust $\gamma$ in response to local curvature information, potentially bridging the gap between the static convergence guarantees established here and practical online parameter tuning.
\bibliographystyle{elsarticle-harv}       
\bibliography{autosam}           



\appendix
\section{Proof Details} 
\subsection{Proof of Lemma~\ref{lem:quad_constraint}}
At equilibrium, we have $\boldsymbol{y}_* = \boldsymbol{x}_*$ and $\boldsymbol{u}_* = \nabla f(\boldsymbol{x}_*) = \mathbf{0}$. Thus, the error vectors are $\tilde{\boldsymbol{y}}_k = \boldsymbol{y}_k - \boldsymbol{x}_*$ and $\tilde{\boldsymbol{u}}_k = \nabla f(\boldsymbol{y}_k)$. Substituting these into the RC~\eqref{eq:RC_ineq} yields
\begin{equation*}
	2\langle \tilde{\boldsymbol{y}}_k, \tilde{\boldsymbol{u}}_k \rangle \ge \mu\|\tilde{\boldsymbol{u}}_k\|^2 + \lambda\|\tilde{\boldsymbol{y}}_k\|^2,
\end{equation*}
which directly rearranges to the quadratic form~\eqref{eq:wrc}.
\subsection{Proof of Lemma~\ref{lem:stability}}   
Suppose the LMI~\eqref{eq:LMI} holds for $\mathbf{P} \succ 0$ and $\rho \in (0,1)$. Define the Lyapunov function
\begin{equation*}
	V(\boldsymbol{\phi}_k) = (\boldsymbol{\phi}_k - \boldsymbol{\phi}_*)^\top (\mathbf{P} \otimes \mathbf{I}_d) (\boldsymbol{\phi}_k - \boldsymbol{\phi}_*).
\end{equation*}
Multiplying LMI~\eqref{eq:LMI} by $[\,(\boldsymbol{\phi}_k - \boldsymbol{\phi}_*)^\top, (\boldsymbol{u}_k - \boldsymbol{u}_*)^\top\,]^\top$ and its transpose yields
\begin{equation}\label{eq:lmi_expanded}
	\begin{aligned}
		&\begin{bmatrix} \boldsymbol{\phi}_k - \boldsymbol{\phi}_* \\ \boldsymbol{u}_k - \boldsymbol{u}_* \end{bmatrix}^\top
		\left( \begin{bmatrix}
			\mathbf{A}^\top \mathbf{P} \mathbf{A}-\rho^2 \mathbf{P} & \mathbf{A}^\top \mathbf{P} \mathbf{B} \\
			\mathbf{B}^\top \mathbf{P} \mathbf{A} & \mathbf{B}^\top \mathbf{P} \mathbf{B}
		\end{bmatrix} \otimes \mathbf{I}_d \right)
		\begin{bmatrix} \boldsymbol{\phi}_k - \boldsymbol{\phi}_* \\ \boldsymbol{u}_k - \boldsymbol{u}_* \end{bmatrix} \\
		&\quad + 
		\begin{bmatrix} \boldsymbol{y}_k - \boldsymbol{y}_* \\ \boldsymbol{u}_k - \boldsymbol{u}_* \end{bmatrix}^\top
		(\mathbf{M} \otimes \mathbf{I}_d)
		\begin{bmatrix} \boldsymbol{y}_k - \boldsymbol{y}_* \\ \boldsymbol{u}_k - \boldsymbol{u}_* \end{bmatrix} \le 0,
	\end{aligned}
\end{equation}
where $\boldsymbol{y}_k - \boldsymbol{y}_* = (\mathbf{C}\otimes \mathbf{I}_d)(\boldsymbol{\phi}_k - \boldsymbol{\phi}_*) + (\mathbf{D}\otimes \mathbf{I}_d)(\boldsymbol{u}_k - \boldsymbol{u}_*)$. By Lemma~\ref{lem:quad_constraint}, the second term satisfies
\begin{equation}\label{eq:wrc_recall}
	\begin{bmatrix} \boldsymbol{y}_k - \boldsymbol{y}_* \\ \boldsymbol{u}_k - \boldsymbol{u}_* \end{bmatrix}^\top
	(\mathbf{M} \otimes \mathbf{I}_d)
	\begin{bmatrix} \boldsymbol{y}_k - \boldsymbol{y}_* \\ \boldsymbol{u}_k - \boldsymbol{u}_* \end{bmatrix} \ge 0.
\end{equation}
Thus, the first term in~\eqref{eq:lmi_expanded} is non-positive. Substituting the state transition $\boldsymbol{\phi}_{k+1} - \boldsymbol{\phi}_* = (\mathbf{A}\otimes \mathbf{I}_d)(\boldsymbol{\phi}_k - \boldsymbol{\phi}_*) + (\mathbf{B}\otimes \mathbf{I}_d)(\boldsymbol{u}_k - \boldsymbol{u}_*)$ into its expansion yields
\begin{equation*}
	V(\boldsymbol{\phi}_{k+1}) - \rho^2 V(\boldsymbol{\phi}_k) \le 0.
\end{equation*}
Recursively applying this inequality gives $V(\boldsymbol{\phi}_{k}) \le \rho^{2k} V(\boldsymbol{\phi}_0)$. By the Rayleigh quotient for $\mathbf{P} \otimes \mathbf{I}_d$, we have
\begin{equation}\label{eq:Ray}
	\begin{aligned}
		\lambda_{\min}(\mathbf{P})\|\boldsymbol{\phi}_k - \boldsymbol{\phi}_*\|^2 &\le V(\boldsymbol{\phi}_k) \le \rho^{2k} V(\boldsymbol{\phi}_0) \\
		&\le \rho^{2k}\lambda_{\max}(\mathbf{P})\|\boldsymbol{\phi}_0 - \boldsymbol{\phi}_*\|^2.
	\end{aligned}
\end{equation}
Dividing by $\lambda_{\min}(\mathbf{P})$ and taking the square root yields
\begin{equation*}
	\|\boldsymbol{\phi}_k - \boldsymbol{\phi}_*\| \le \sqrt{\kappa(\mathbf{P})} \rho^k \|\boldsymbol{\phi}_0 - \boldsymbol{\phi}_*\|.
\end{equation*}
\subsection{Proof of Lemma~\ref{lem:relaxed_qc}}
From $\boldsymbol{v}_k = -\gamma\nabla f(\boldsymbol{y}_k) - \delta \boldsymbol{y}_{k}$, we have $\nabla f(\boldsymbol{y}_k) = -\frac{1}{\gamma}(\boldsymbol{v}_k + \delta \boldsymbol{y}_k)$. Since $\nabla f(\boldsymbol{y}_*) = \mathbf{0}$, the equilibrium satisfies $\boldsymbol{v}_* = -\delta \boldsymbol{y}_*$. Using the error variables $\tilde{\boldsymbol{y}}_k = \boldsymbol{y}_k - \boldsymbol{y}_*$ and $\tilde{\boldsymbol{v}}_k = \boldsymbol{v}_k - \boldsymbol{v}_*$, the gradient is
\begin{equation*}
	\nabla f(\boldsymbol{y}_k) = -\frac{1}{\gamma}(\boldsymbol{v}_k + \delta \boldsymbol{y}_k) + \frac{1}{\gamma}(\boldsymbol{v}_* + \delta \boldsymbol{y}_*) = -\frac{1}{\gamma}(\tilde{\boldsymbol{v}}_k + \delta \tilde{\boldsymbol{y}}_k).
\end{equation*}
Substituting this into the RC~\eqref{eq:RC_ineq} yields
\begin{equation*}
	\begin{aligned}
		-\frac{2}{\gamma}\langle \tilde{\boldsymbol{y}}_k, \tilde{\boldsymbol{v}}_k \rangle - \frac{2\delta}{\gamma}\|\tilde{\boldsymbol{y}}_k\|^2 &\ge \mu\left\| \frac{1}{\gamma}(\tilde{\boldsymbol{v}}_k + \delta \tilde{\boldsymbol{y}}_k) \right\|^2 + \lambda\|\tilde{\boldsymbol{y}}_k\|^2 \\
		&= \frac{\mu}{\gamma^2}\|\tilde{\boldsymbol{v}}_k\|^2 + \frac{2\mu\delta}{\gamma^2}\langle \tilde{\boldsymbol{y}}_k, \tilde{\boldsymbol{v}}_k \rangle \\
		&\quad + \left(\frac{\mu\delta^2}{\gamma^2} + \lambda\right)\|\tilde{\boldsymbol{y}}_k\|^2.
	\end{aligned}
\end{equation*}
Multiplying by $\gamma^2$ and rearranging gives
\begin{equation*}
	(-\lambda\gamma^2 - 2\gamma\delta - \mu\delta^2)\|\tilde{\boldsymbol{y}}_k\|^2 + 2(-\gamma - \mu\delta)\langle \tilde{\boldsymbol{y}}_k, \tilde{\boldsymbol{v}}_k \rangle - \mu\|\tilde{\boldsymbol{v}}_k\|^2 \ge 0,
\end{equation*}
which gives the quadratic form~\eqref{eq:relaxed_qc}.
\subsection{Proof of Theorem~\ref{thm:param_range}}

We verify the three prerequisites of $\text{KYPC}(\mathbf{\hat{A}}, \mathbf{\hat{M}})$ established in Corollary~\ref{coro:kypc}.

\textit{Prerequisite 1.} Requires $\det(e^{j\omega}\mathbf{I} - \mathbf{\hat{A}}) \neq 0$ for all $\omega \in \mathbb{R}$, meaning $\mathbf{\hat{A}}$ has no eigenvalues on the unit circle. This is trivially satisfied by the strict Schur stability guaranteed in Prerequisite 2.

\textit{Prerequisite 2.} Requires $\mathbf{\hat{A}}$ to be Schur stable. Since $\mathbf{\hat{A}}$ defined in~\eqref{eq:relaxed_matrices} has a trivial eigenvalue at $z=0$, its stability depends entirely on its leading $2 \times 2$ principal submatrix $\mathbf{\hat{A}}_{2\times2}$. Its characteristic polynomial is $p(z) = z^2 - \text{tr}(\mathbf{\hat{A}}_{2\times2})z + \det(\mathbf{\hat{A}}_{2\times2})$, where $\text{tr}(\mathbf{\hat{A}}_{2\times2}) = 1 + \beta + \frac{(\alpha+\gamma)\delta}{\gamma}$ and $\det(\mathbf{\hat{A}}_{2\times2}) = \beta + \delta$. By the Jury stability criterion, $\mathbf{\hat{A}}_{2\times2}$ is Schur stable if and only if:
\begin{enumerate}
	\item $p(1) > 0 \implies -\frac{\alpha \delta}{\gamma} > 0$. Given $\alpha, \gamma > 0$, this explicitly requires $\delta < 0$.
	\item $p(-1) > 0 \implies 2(1+\beta) + \delta\left(\frac{2\gamma+\alpha}{\gamma}\right) > 0$.
	\item $|\det(\mathbf{\hat{A}}_{2\times2})| < 1 \implies -1 - \beta < \delta < 1 - \beta$. Intersecting this valid negative interval with the condition $\delta < 0$ dictates $\beta \in (-1, 1)$.
\end{enumerate}
Combining these constraints establishes the feasible interval for the relaxation parameter
\begin{equation*}
	\delta \in \left(\max\left\{-1-\beta, -\frac{2\gamma(1+\beta)}{2\gamma+\alpha}\right\}, 0\right).
\end{equation*}

\textit{Prerequisite 3.} Requires the upper-left block $r_1$ of $\mathbf{\hat{M}}$ to be positive semidefinite. Enforcing $r_1 \ge 0$ yields $-(\lambda \gamma^2 + 2\delta \gamma +\mu\delta^2) \ge 0$. Since $\delta < 0$, solving this quadratic inequality for $\gamma$ provides the bound
\begin{equation}\label{eq:root_bound}
	\frac{-\delta(1-\sqrt{1-\lambda\mu})}{\lambda}\leq \gamma \leq\frac{-\delta(1+\sqrt{1-\lambda\mu})}{\lambda}.
\end{equation}

A valid relaxation parameter $\delta < 0$ exists if and only if the feasible intervals from Prerequisites 2 and 3 intersect. Rearranging~\eqref{eq:root_bound} yields the lower bound $-\delta \ge \frac{\lambda \gamma}{1+\sqrt{1-\mu\lambda}}$. Simultaneously, the binding constraint from Prerequisite 2 dictates the upper bound $-\delta < \frac{2\gamma(1+\beta)}{2\gamma+\alpha}$. Enforcing a strictly non-empty intersection between these bounds gives
\begin{equation*}
	\frac{\lambda \gamma}{1+\sqrt{1-\mu\lambda}} < \frac{2\gamma(1+\beta)}{2\gamma+\alpha}.
\end{equation*}
Canceling $\gamma > 0$ and rearranging algebraically simplifies the inequality to $2\gamma + \alpha < \frac{2(1+\beta)(1+\sqrt{1-\mu\lambda})}{\lambda}$. Constraining $\alpha, \gamma > 0$ explicitly produces the parameter boundaries presented in~\eqref{eq:param_condition}.
\subsection{Proof of Lemma~\ref{lem:local_rc_init}}

For $\boldsymbol{\phi}_0 = [\boldsymbol{x}_0^\top, \boldsymbol{x}_{-1}^\top, \boldsymbol{v}_{-1}^\top]^\top$ and $\boldsymbol{\phi}_* = [\boldsymbol{x}_*^\top, \boldsymbol{x}_*^\top, \boldsymbol{v}_*^\top]^\top$, the initial deviation expands as
\begin{equation*}
	\|\boldsymbol{\phi}_0 - \boldsymbol{\phi}_*\|^2 = \|\boldsymbol{x}_0 - \boldsymbol{x}_*\|^2 + \|\boldsymbol{x}_{-1} - \boldsymbol{x}_*\|^2 + \|\boldsymbol{v}_{-1} - \boldsymbol{v}_*\|^2.
\end{equation*}
Since $\boldsymbol{v}_{-1} = -\gamma \nabla f(\boldsymbol{x}_{-1}) - \delta \boldsymbol{x}_{-1}$ and $\nabla f(\boldsymbol{x}_*) = \mathbf{0}$, the velocity error satisfies
\begin{equation*}
	\begin{aligned}
		\|\boldsymbol{v}_{-1} - \boldsymbol{v}_*\| &= \|-\gamma (\nabla f(\boldsymbol{x}_{-1}) - \nabla f(\boldsymbol{x}_*)) - \delta (\boldsymbol{x}_{-1} - \boldsymbol{x}_*)\| \\
		&\le \gamma \|\nabla f(\boldsymbol{x}_{-1})\| + |\delta| \|\boldsymbol{x}_{-1} - \boldsymbol{x}_*\|.
	\end{aligned}
\end{equation*}
Applying the Cauchy-Schwarz inequality to the local RC yields $\|\nabla f(\boldsymbol{x}_{-1})\| \le \frac{2}{\mu} \|\boldsymbol{x}_{-1} - \boldsymbol{x}_*\|$. Substituting this bound limits the velocity error to
\begin{equation*}
	\|\boldsymbol{v}_{-1} - \boldsymbol{v}_*\| \le \left(\frac{2\gamma}{\mu} + |\delta|\right) \|\boldsymbol{x}_{-1} - \boldsymbol{x}_*\|.
\end{equation*}
Given $\boldsymbol{x}_{-1}, \boldsymbol{x}_0 \in \mathcal{B}_{r_0}(\boldsymbol{x}_*)$, the initial distance of the augmented state is bounded by
\begin{equation*}
	\|\boldsymbol{\phi}_0 - \boldsymbol{\phi}_*\| \le r_0 \sqrt{2 + \left(\frac{2\gamma}{\mu} + |\delta|\right)^2} = \frac{\epsilon \|\boldsymbol{x}_*\|}{\sqrt{\kappa(\mathbf{P})}}.
\end{equation*}
Invoking the state error bound~\eqref{eq:Ray} established in Lemma~\ref{lem:stability}, and noting $\boldsymbol{y}_k = \mathbf{\hat{C}}\boldsymbol{\phi}_k$ with $\mathbf{\hat{C}} = [\mathbf{I}_d, \mathbf{0}, \mathbf{0}]$, the output error satisfies
\begin{equation*}
	\|\boldsymbol{y}_k - \boldsymbol{x}_*\| \le \|\boldsymbol{\phi}_k - \boldsymbol{\phi}_*\| \le \sqrt{\kappa(\mathbf{P})} \|\boldsymbol{\phi}_0 - \boldsymbol{\phi}_*\| \le \epsilon \|\boldsymbol{x}_*\|.
\end{equation*}
This confirms that the trajectory resides entirely within the local RC domain.
\subsection{Proof of Lemma~\ref{lem:transfer_function}}
Given the matrices in~\eqref{eq:relaxed_matrices} and~\eqref{eq:sys_scalars}, and $\mathbf{\hat{D}} = \mathbf{0}$, the transfer function is $H(z) = \mathbf{\hat{C}}(z\mathbf{I} - \mathbf{\hat{A}})^{-1}\mathbf{\hat{B}}$. Let $\mathbf{v} = (z\mathbf{I} - \mathbf{\hat{A}})^{-1}\mathbf{\hat{B}} = [v_1, v_2, v_3]^\top$. Since $\mathbf{\hat{C}} = [1, 0, 0]$, we have $H(z) = v_1$. The vector $\mathbf{v}$ satisfies $(z\mathbf{I} - \mathbf{\hat{A}})\mathbf{v} = \mathbf{\hat{B}}$, which expands to
\begin{equation*}
	\left\{
	\begin{aligned}
		(z - A_{11})v_1 - A_{12}v_2 + v_3 &= B_1, \\
		-v_1 + zv_2 &= 0, \\
		zv_3 &= 1.
	\end{aligned}
	\right.
\end{equation*}
Solving bottom-up yields $v_3 = z^{-1}$ and $v_2 = z^{-1}v_1$. Substituting these into the first equation gives
\begin{equation*}
	(z - A_{11})v_1 - A_{12}(z^{-1}v_1) + z^{-1} = B_1.
\end{equation*}
Multiplying by $z$ and isolating $v_1$ directly yields the scalar transfer function~\eqref{eq:transfer_func_explicit}.
\subsection{Proof of Theorem~\ref{thm:fdi_polynomial}}
By the KYP lemma, the LMI~\eqref{eq:LMI_generic} is equivalent to the following FDI for all $\omega \in \mathbb{R}$:
\begin{equation}\label{eq:fdi_general}
	\left[ \begin{smallmatrix} (e^{j \omega} \mathbf{I}-\mathbf{\hat{A}})^{-1} \mathbf{\hat{B}} \\ 1 \end{smallmatrix} \right]^*
	\mathbf{\Pi}
	\left[ \begin{smallmatrix} (e^{j \omega} \mathbf{I}-\mathbf{\hat{A}})^{-1} \mathbf{\hat{B}} \\ 1 \end{smallmatrix} \right] < 0,
\end{equation}
where $\mathbf{\Pi} = \left[\begin{smallmatrix} \mathbf{\hat{C}}^\top & \mathbf{0} \\ \mathbf{0} & 1 \end{smallmatrix}\right] \mathbf{\hat{M}} \left[\begin{smallmatrix} \mathbf{\hat{C}} & \mathbf{0} \\ \mathbf{0} & 1 \end{smallmatrix}\right]$. Noting $\mathbf{\hat{C}}(e^{j \omega} \mathbf{I}-\mathbf{\hat{A}})^{-1} \mathbf{\hat{B}} = H(e^{j\omega})$, the FDI reduces to the scalar inequality $r_1 |H(z)|^2 + 2 r_2 \text{Re}\big(H(z)\big) - \mu < 0$, where $z = e^{j\omega}$.
By Lemma~\ref{lem:transfer_function}, $H(z) = N(z)/D(z)$ with $N(z) = B_1 z - 1$ and $D(z) = z^2 - A_{11}z - A_{12}$. Since Schur stability ensures $|D(z)|^2 > 0$ on the unit circle, multiplying the scalar inequality by $|D(z)|^2$ clears the denominator:
\begin{equation}\label{eq:cleared_fdi}
	r_1 |N(z)|^2 + 2r_2 \text{Re}\big(N(z) D^*(z)\big) - \mu |D(z)|^2 < 0.
\end{equation}
On the unit circle $z^* = z^{-1}$, the conjugate products expand as $|N(z)|^2 = N(z)N(z^{-1})$ and $2\text{Re}(N(z)D^*(z)) = N(z)D(z^{-1}) + N(z^{-1})D(z)$. Substituting these relations alongside $z+z^{-1} = 2x$ and $z^2+z^{-2} = 4x^2 - 2$ (where $x = \cos\omega \in [-1, 1]$) into~\eqref{eq:cleared_fdi}, the inequality algebraically groups into the quadratic polynomial:
\begin{equation*}
	\Gamma_2 x^2 + \Gamma_1 x + \Gamma_0 < 0.
\end{equation*}
By substituting the definitions in~\eqref{eq:sys_scalars} into the grouped coefficients, the algebraic simplification explicitly cancels out the relaxation parameter $\delta$, directly yielding the expressions in~\eqref{eq:gamma_coeffs}.
\subsection{Proof of Theorem~\ref{thm:analytical_convergence_region}}

By Theorem~\ref{thm:fdi_polynomial}, the condition $P(x) < 0$ on $[-1, 1]$ is analyzed via the sign of $\Gamma_2$.

\noindent\textbf{Case 1: $\Gamma_2 \ge 0$.} This requires $\gamma \ge \mu\beta$. As $P(x)$ is convex, $P(x) < 0$ on $[-1, 1]$ holds if and only if $P(1) < 0$ and $P(-1) < 0$. Given $\alpha, \lambda > 0$, the right endpoint $P(1) = -\lambda\alpha^2 < 0$ holds trivially. For the left endpoint, letting $S = \alpha + 2\gamma$ yields
\begin{equation*}
	P(-1) = -\lambda S^2 + 4(1+\beta)S - 4\mu(1+\beta)^2 < 0.
\end{equation*}
Solving this quadratic inequality for $S$ yields $S < S_1$ or $S > S_2$, where
\begin{equation*}
	S_{1,2} = \frac{2(1+\beta)(1 \mp \sqrt{1-\mu\lambda})}{\lambda}.
\end{equation*}
Since Theorem~\ref{thm:param_range} explicitly constrains $S < S_2$, the feasible interval is truncated to the left branch:
\begin{equation}\label{eq:x12}
	\alpha + 2\gamma < \frac{2(1+\beta)(1-\sqrt{1-\mu\lambda})}{\lambda}.
\end{equation}
For a non-empty feasible interval $\alpha \in(0, S_1-2gamma)$, we require $2\gamma < S_1$. Rationalizing the numerator elegantly simplifies this bound to
\begin{equation*}
	\gamma < \frac{(1+\beta)(1-\sqrt{1-\mu\lambda})}{\lambda} = \frac{\mu(1+\beta)}{1+\sqrt{1-\mu\lambda}}.
\end{equation*}
Intersecting this upper bound with $\gamma \ge \mu\beta$ yields the over-damped region $\Omega_u$ in~\eqref{eq:region_u}.

\noindent\textbf{Case 2: $\Gamma_2 < 0$.} This requires $0 < \gamma < \mu\beta$, implying $\beta \in (0, 1)$. The necessary endpoint conditions $P(\pm 1) < 0$ again yield the strict bound~\eqref{eq:x12}. From the coefficients in~\eqref{eq:gamma_coeffs}, the axis of symmetry is
\begin{equation*}
	x_v = -\frac{\lambda\gamma(\alpha+\gamma) - (1+\beta)(\alpha + 2\gamma) + \mu(1+\beta)^2}{4(\gamma - \mu\beta)}.
\end{equation*}
When $x_v \in [-1, 1]$, the global maximum must satisfy $P(x_v) < 0$, giving
\begin{equation*}
	\begin{aligned}
		P(x_v) &= -\frac{\Gamma_1^2 - 4\Gamma_2\Gamma_0}{4\Gamma_2} < 0 \\
		&\Leftrightarrow \Delta(\alpha) := \Gamma_1^2 - 4\Gamma_2\Gamma_0 = c_2 \alpha^2 + c_1 \alpha + c_0 < 0,
	\end{aligned}
\end{equation*}
where the explicit polynomial coefficients of $\alpha$ are derived as
\begin{equation*}
	\begin{aligned}
		c_2 &= 4\lambda^2\gamma^2 + 8\lambda(1-\beta)\gamma + 4[(1+\beta)^2 - 4\mu\lambda\beta], \\
		c_1 &= 8\lambda^2\gamma^3 - 8[3\lambda(1+\beta) + 4]\gamma^2 \\
		&\quad + 8[\lambda\mu(1+\beta)^2 + 4\mu\beta - 2(1-\beta^2)]\gamma \\
		&\quad - 8\mu(1+\beta)(1-\beta)^2, \\
		c_0 &= 4\lambda^2\gamma^4 - 16[\lambda(1+\beta) + 2]\gamma^3 \\
		&\quad + 8[2(1-\beta)^2 + \lambda\mu(1+\beta)^2 + 4\mu\beta]\gamma^2 \\
		&\quad + 16\mu(1+\beta+3\beta^2-\beta^3)\gamma \\
		&\quad + 4\mu^2(1+\beta)^2(1-6\beta+\beta^2).
	\end{aligned}
\end{equation*}
This defines the implicit upper bound $\Xi_d(\beta, \gamma) := \{ \alpha > 0 \mid \Delta(\alpha) < 0 \}$. Intersecting the universal endpoint restriction~\eqref{eq:x12} with $\Xi_d$ establishes the under-damped region $\Omega_d$ in~\eqref{eq:region_d}. The overall feasible region is the union $\Omega = \Omega_u \cup \Omega_d$.

\end{document}